\documentclass{elsart}
\usepackage{natbib}
\usepackage{graphics}
\usepackage{epsfig}
\usepackage{amsfonts}
\usepackage{amssymb}
\usepackage{amsmath}

\newcommand{\euO}{\mathfrak O}
\newcommand{\euP}{\mathfrak P}

\newcommand{\bQ}{\mathbb Q}
\newcommand{\bZ}{\mathbb Z}
\newcommand{\bF}{\mathbb F}

\begin{document}

\begin{frontmatter}

\title{On the necessity of new ramification breaks}

\author[label1]{Nigel P. Byott}, \ead{N.P.Byott@ex.ac.uk}
\author[label2]{G. Griffith Elder\corauthref{cor1}\thanksref{label3}}
\ead{elder@vt.edu} \corauth[cor1]{Corresponding Author}
\address[label1]{School of Engineering, Computer Science and
Mathematics, University of Exeter, Exeter EX4 4QE U.K.}
\address[label2]{Department of Mathematics, Virginia Tech, Blacksburg
VA 24061-0123 U.S.A.}  \thanks[label3]{Elder was partially supported
by the National Science Foundation under Grant No. 0201080.}

\begin{abstract}
Ramification invariants are necessary, but not in general sufficient,
to determine the Galois module structure of ideals in local number
field extensions. This insufficiency is associated with elementary
abelian extensions, where one can define a refined ramification
filtration -- one with more ramification breaks
\cite{elder:newbreaks}. The first refined break number comes from the
usual ramification filtration and is therefore necessary. Here we
study the second refined break number.
\end{abstract}

\bibliographystyle{amsalpha} 

\end{frontmatter}
\section{Introduction}

Let $p$ be a prime integer, and let $K$ be a finite extension of the
field $\bQ_p$ of $p$-adic numbers, with absolute ramification index
$e_K$ and inertia degree $f$. Let $N$ be a finite, fully ramified,
Galois $p$-extension of $K$, let $G=\mbox{Gal}(N/K)$, and let $\euP_N$
be the maximal ideal of the valuation ring $\euO_N$ of $N$.  Also, let
$T$ be the maximal unramified subfield of $K$. Thus the valuation ring
$\euO_T$ of $T$ is the ring of Witt vectors of $\mathbb{F}_q$, where
$q=p^f$. It is natural to ask about the structure of each ideal
$\euP_N^r$ under the canonical action of the group ring $\euO_T[G]$.
This question has its roots in the Normal Basis Theorem, see e.g.
\cite[p.~344]{lang}, and in the Normal Integral Basis Theorem of
E. Noether \cite{emmy}.

Complexity, however, threatens to overwhelm any complete, explicit
description, even when one restricts oneself to relatively simple
Galois groups \cite{elder:annals, elder:bord, elder:c8}. So instead, we ask
for those invariants upon which the structure depends.  Certainly
these must include those associated with the usual ramification
filtration
\[G_i=\{\sigma\in G: (\sigma-1)\euP_N\subset\euP_N^{i+1}\}.\]
For example, it is easily shown that the ramification break numbers
(that is, the integers $b$ such that $G_b\supsetneq G_{b+1}$) are
necessary to determine the Galois module structure of the ideals of
$\euO_N$. To see this, simply consider the structure of the ideal
fixed by $G_{b+1}$, namely $(\euP_N^r)^{G_{b+1}}$, over the group ring
$\euO_T[\sigma]$ for some $\sigma\in G_b\setminus G_{b+1}$. Since
$G_b/G_{b+1}\cong C_p^s$ is necessarily elementary abelian \cite[IV
\S2 Prop 7 Cor 3]{serre:local}, this is a module over the cyclic group
ring $\euO_T[C_p]$, for which there are exactly three indecomposable
modules: the trivial module $\euO_T$, the group ring or regular
representation $\euO_T[C_p]$, and the module $\euO_T[\zeta_p]$ where
$\sigma$ acts via multiplication by the $p$-th root of unity $\zeta_p$
\cite[Thm 34.31]{reiner}. Now proceed as in \cite[Thm 1]{martha} to
see how the multiplicities of these three modules are parametrized by
$b$ (along with the absolute ramification degree).

The usual ramification invariants are not however sufficient to
determine the Galois module structure of ideals. This was observed in
\cite{elder:onebreak} where we considered biquadratic extensions (the
case $p=2$) with one break. The work presented here, together with
\cite{elder:newbreaks}, stems from our ongoing effort to fully
understand the implications of that paper, and to extend its results
to arbitrary $p$.  With hindsight we can now say that the
insufficiency of the usual ramification filtration is tied to the
elementary abelian quotients of consecutive ramification groups
$G_b/G_{b+1}$, but that there is a `repair'. We can focus on the
elementary abelian extension with Galois group $G_b/G_{b+1}$ and
define a new refined ramification filtration, one with more
information -- more breaks \cite{elder:newbreaks}. In this paper, we
amend the definition from \cite{elder:newbreaks} slightly; study the
necessity, for the Galois module structure of ideals, of the first
piece of new information that this refined ramification filtration
provides -- the second refined break; and explicitly describe the
Galois module structure of ideals in bicyclic extensions under {\em
maximal refined ramification}, when this second refined break achieves
a natural upper bound.

\subsection{Refined Ramification Filtration}

Let $N/K$ be a fully ramified, elementary abelian $p$-extension with
one break in its ramification filtration, at $b$. So
$G=\mbox{Gal}(N/K)\cong G_b/G_{b+1}$.  Note that $G$ is a vector space
over $\bF_p$, the field with $p$ elements.  To enable the residue
field $\bF_q$ to act on $G$ as well, let $\bZ_{(p)}$ denote the integers
localized at $p$, and define {\em truncated exponentiation} by the
polynomial
\[(1+X)^{[Y]}:=\sum_{i=0}^{p-1}\binom{Y}{i}X^i\in \bZ_{(p)}[X,Y],\] a
truncation of the usual binomial series.  Now let
$\mathcal{A}=(\sigma-1:\sigma\in G)$ denote the augmentation ideal of
$\euO_T[G]$.  If $L$ is any finite extension of $T$ that is contained
in $K$, then $\euO_L\mathcal{A}$ is the augmentation ideal of
$\euO_L[G]$.  For any $\kappa\in \euO_L$ and any $x\in
1+\euO_L\mathcal{A}$, truncated exponentiation gives a well-defined
element $x^{[\kappa]}$ of $ 1+\euO_L\mathcal{A}$. This does not make $
1+\euO_L\mathcal{A}$ into an $\euO_L$-module since, for example, we do
not in general have $(x^{[\kappa]})^{[\kappa']}=x^{[\kappa\kappa']}$.
To address this problem we could choose to work with the
quotient group $(1+\euO_L\mathcal{A})/(1+p\euO_L\mathcal{A})$.

This is the approach of \cite{elder:newbreaks} in the case $L=T$,
where we proposed working with the quotient group
$(1+\mathcal{A})/(1+p\mathcal{A})$ over the field
$\euO_T/p\euO_T=\mathbb{F}_q$. As noted there,
$(1+\mathcal{A})/(1+p\mathcal{A})$ is a ``near-space'' over
$\mathbb{F}_q$: it satisfies all the properties of a vector space over
$\mathbb{F}_q$ except the distributive property,
$(x_1x_2)^{[\omega]}\neq x_1^{[\omega]}x_2^{[\omega]}$.  In the case
of biquadratic extensions, the refined ramification filtration of this
near space contains extraneous information in the form of an ``extra''
third refined break \cite[\S4]{elder:newbreaks}. This is undesirable
and expected more generally.

So, in this paper, we propose working with the smaller group
$\mathcal{G}=(1+\mathcal{A})/(1+\mathcal{A}^p)$. Notice that because
$G$ is elementary abelian, we have
$p\mathcal{A}\subset\mathcal{A}^p$. Following \cite[Thm
2.1]{elder:newbreaks} and the discussion leading to \cite[Cor
2.3]{elder:newbreaks}, we find
\[(\omega,x)\in \bF_q\times \mathcal{G}\longrightarrow x^{[\omega]}\in
\mathcal{G}\] 
is an $\bF_q$-action that
endows $\mathcal{G}$ with the structure of an
$\bF_q$-vector space. 
Let $G^{\bF}$ be the span of the image of $G$
in $\mathcal{G}$. Clearly
\[G^{\bF}\cong\bF_q\otimes_{\bF_p}G.\]

Now choose any $\alpha\in N$ with $v_N(\alpha)=b$.  Because of
\cite[Cor 4]{elder:valNB}, such elements generate normal field bases
and are thus valuable for Galois module structure.  Following the
treatment of the usual ramification filtration
\cite[p62]{serre:local}, define a function $i_\alpha$ on $\bar{x}\in
G^{\bF}$ by the formula $i_\alpha(\bar{x})=\sup\{v_N((x-1)\alpha):x\in
1+\mathcal{A},\; x\cdot(1+\mathcal{A}^p)=\bar{x}\}$. The refined
ramification filtration of $G^{\bF}$, which cannot as yet be
considered canonical as it apparently depends upon a choice of
$\alpha$, is defined by
\[G_j^{\mathcal{F},\alpha}=\{\bar{x}\in G^{\bF}:
i_{\alpha}(\bar{x})\geq v_N(\alpha)+j\}.\] This leads to a definition
of {\em refined breaks}: integers $j$ such that $G^{\mathcal{F},
\alpha}_j\supsetneq G^{\mathcal{F}, \alpha}_{j+1}$.  Because of
\cite[Cor 4]{elder:valNB} and by following \cite[Thm
3.3]{elder:newbreaks}, we see that there are exactly $\log_p|G|$
refined breaks.

The value of the first refined break is $b$ (the usual ramification
number) and so is clearly necessary for Galois module structure.  The
purpose of this paper is threefold:
\begin{enumerate}
\item Show that the second refined break, which we call $b_*$, is canonical.
\item Characterize those integers that appear as $b_*$ in some
extension.
\item Discuss the
relevance of $b_*$ for Galois module structure.
\end{enumerate}

Notice that we can repeat the procedure that was just described for
each bicyclic subgroup $H\cong C_p^2$ of $G$. In each case there will
be two refined breaks: $b$ and a second refined break $b_H$.  Since
the second refined break associated with $G$ is the minimum of these
$b_H$, there is a bicyclic subgroup $H$ with the refined breaks
$b<b_*$.  We can restrict our attention to this particular bicyclic
extension and answer all three questions.  Since the implications for
the general Galois extension should be clear, we henceforth restrict
our attention to $N/K$, a bicyclic extension with
$G=\mbox{Gal}(N/K)\cong C_p^2$ and refined breaks $b<b_*$.

\subsection{Outline}
In \S2 we determine the value of $b_*$, find that it is canonical and
moreover, that it satisfies $b<b_*\leq pb$ with the additional
condition that $b_*\equiv b \bmod p$ when $b_*<pb$. The special case
when $b_*=pb$ will be called {\em maximal refined ramification} (MRR)
and $(p-1+1/p)b<b_*<pb$, {\em near maximal refined ramification}
(NMRR).  In \S3 prove two results in Galois module structure. We find
in Theorem 12 of \S3.1 that outside of NMRR, the $\bF_q[G]$-structure
of $\euP_N^r/p\euP_N^r$ {\em depends} upon $b_*$, and therefore so too
does the $\euO_T[G]$-structure of $\euP_N^r$.  This addresses the
question raised in the title of this paper by proving that the second
refined break is necessary for the Galois module structure of ideals,
as long as it is ``not too big'' relative to $b$. Then in \S3.2 we show
in Theorem 18 how MRR allows an easy, rather transparent and explicit
description of Galois module structure in terms of $\euO_T[G]$-ideals.

\section{The Refined Ramification Filtration in Bicyclic Extensions} 
Let $N/K$ be a fully ramified bicyclic extension with
$G=\mbox{Gal}(N/K)\cong C_p^2$ and one ramification break at $b$,
which necessarily satisfies $0<b<pe_K/(p-1)$ and $\gcd(b,p)=1$. We
begin a process now that will define an integer, our candidate for the
second refined break.

Choose $\rho_0\in N$ with $v_N(\rho_0)=b$, and choose a pair of
generators $\gamma, \sigma$, so that $G=\langle\gamma, \sigma\rangle$.
Since $v_N((\gamma-1)\rho_0)=v_N((\sigma-1)\rho_0)=2b$ and $N/T$ is
fully ramified, there is a unique $p^f-1$ root of unity
$\omega_{\gamma,\sigma}$ such that $(\gamma-1)\rho_0\equiv
\omega_{\gamma,\sigma}(\sigma-1)\rho_0 \bmod\euP_N^{2b+1}$.  Since
$\gamma\not\in\langle\sigma\rangle$, $\omega_{\gamma,\sigma}^{p-1}\neq
1$. But how does $\omega_{\gamma,\sigma}$ depend upon our choice of
group generators?  Observe
$(\gamma^i-1)=i(\gamma-1)+\sum_{j=2}^i\binom{i}{j}(\gamma-1)^j$ and
that similarly $(\sigma^j-1)\equiv j(\sigma-1)$ modulo higher powers
of $(\sigma-1)$.  Moreover
$(\gamma^i\sigma^j-1)=(\gamma^i-1)+(\sigma^j-1)+(\gamma^i-1)(\sigma^j-1)$. As
a result, $(\gamma^i\sigma^j-1)\rho_0\equiv i(\gamma-1)\rho_0
+j(\sigma-1)\rho_0 \equiv
(i\omega_{\gamma,\sigma}+j)(\sigma-1)\rho_0\bmod\euP_N^{2b+1}$.  This
means that the change of group generators $\langle\gamma,
\sigma\rangle=\langle\gamma^a\sigma^b,\gamma^c\sigma^d\rangle$,
resulting from
$$\begin{bmatrix}a&b\\c&d\end{bmatrix}\in \mbox{GL}_2(\mathbb{F}_p),$$
leads to 
$(\gamma^a\sigma^b-1)\rho_0\equiv \omega_{\gamma^a\sigma^b,\gamma^c\sigma^d}(\gamma^c\sigma^d-1)\rho_0 \bmod\euP_N^{2b+1}$
and thus
$a\omega_{\gamma,\sigma}+b
\equiv \omega_{\gamma^a\sigma^b,\gamma^c\sigma^d}
(c\omega_{\gamma,\sigma}+d)
\bmod\euP_T$.
In other words, if we identify the $p^f-1$ roots
of unity with the nonzero elements of the finite field $\mathbb{F}_q$,
we have
$$\omega_{\gamma^a\sigma^b,\gamma^c\sigma^d}=
\frac{a\omega_{\gamma,\sigma}+b}{c\omega_{\gamma,\sigma}+d}.$$ A
unified approach requires that we identify these roots of unity with
points on the projective line,
$(\omega_{\gamma^a\sigma^b,\gamma^c\sigma^d},1)=
(a\omega_{\gamma,\sigma}+b,c\omega_{\gamma,\sigma}+d)\in
\mathbf{P}^1(\mathbb{F}_q)$. We conclude that while the particular point
$(\omega_{\gamma,\sigma},1)\in\mathbf{P}^1(\mathbb{F}_q)\setminus
\mathbf{P}^1(\mathbb{F}_p)$ depends upon our choice of group
generators, its orbit, $\mbox{Orb}_{N/K}
\subseteq\mathbf{P}^1(\mathbb{F}_q)\setminus
\mathbf{P}^1(\mathbb{F}_p)
$, under
$\mbox{PGL}_2(\mathbb{F}_p)$ is independent of both our choice of
group generators and element $\rho_0$, and should be considered a
basic invariant of the extension.

Fix $\rho_0\in N$ now with $v_N(\rho_0)=b$, and fix our group
generators, so $G=\langle\gamma,\sigma\rangle$.  Rewrite the equation
$(\gamma-1)\rho_0\equiv \omega_{\gamma,\sigma}(\sigma-1)\rho_0
\bmod\euP_N^{2b+1}$ as $\gamma\rho_0\equiv \left
(1+\omega_{\gamma,\sigma}(\sigma-1)\right )\rho_0 \bmod\euP_N^{2b+1}$.
Motivated by the appearance of the first two terms in truncated
exponentiation, we drop subscripts, write
$\omega=-\omega_{\gamma,\sigma}$, and define
$\Theta=\gamma\sigma^{[\omega]}\in\euO_T[G]$. 

Observe that $(\Theta-1)\rho_0\equiv 0\bmod\euP_N^{2b+1}$. Define our
``candidate'' second refined break by
$$b_*:=v_N((\Theta-1)\rho_0)- v_N(\rho_0).$$ This is an integer $> b$,
which may depend upon our choices: of group generators and of
$\rho_0$.  Let $L=N^{\sigma}$ be the fixed field of
$\langle\sigma\rangle$.

The purpose of this paper, as stated in \S1.1, is to address three
goals.  In \S2.1, we address the first goal by proving that $b_*$ is
the second refined break and that it is also canonical (independent of
our choice of $\rho_0$ and also of our choice of the generators for
$G$). In \S2.2, we address the second goal by determining all
realizable second refined breaks.  The third goal is addressed in \S3.

\subsection{The second refined break is canonical}

We begin by establishing the upper bound $b_*\leq pb$.  Recall the
augmentation ideal $\mathcal{A}=(\sigma-1,\gamma-1)\subseteq
\euO_T[G]$ as defined in \S1.1.

\begin{lem} 
Given $a\not\equiv -1\bmod p$, $\rho\in N$ with
$v_{N}(\rho)=(1+ap)b$, $\kappa\in \euO_L$ and $\mu\in
\mathcal{A}^p$. Then $b\leq
v_{N}((\gamma\sigma^{[\kappa]}(1+\mu)-1)\rho)-v_N(\rho)\leq pb$.
\end{lem}

\begin{pf}
We need to prove two inequalities. The first is obvious. So consider
the second inequality and the effect of the trace
$\mbox{Tr}_{N/L}=\Phi_p(\sigma)$ on $\rho$ and on
$\rho_*=(\gamma\sigma^{[\kappa]}(1+\mu)-1)\rho$.  
Because of
\cite[V\S3 Lem 4]{serre:local}, if $v_{N}(\rho_*)>v_N(\rho)+pb$, then
$v_{L}(\mbox{Tr}_{N/L}\rho_*)>v_L(\mbox{Tr}_{N/L}\rho)+b$.

So we prove
$v_{L}(\mbox{Tr}_{N/L}\rho_*)=v_L(\mbox{Tr}_{N/L}\rho)+b$.  Since
$v_N((\sigma-1)\alpha)=v_N(\alpha)+b$ if $\gcd(v_N(\alpha),p)=1$, we
have $v_N((\sigma-1)^{p-1}\rho)=v_N(\rho)+(p-1)b=(1+a)pb\not\equiv 0\bmod p^2$.  Since the
cyclotomic polynomial $\Phi_p(\sigma)\equiv (\sigma-1)^{p-1}\bmod p$
and $(p-1)b<v_N(p)$, we therefore also have
$v_N(\mbox{Tr}_{N/L}\rho)=(1+a)pb$.  So
$\gcd(v_L(\mbox{Tr}_{N/L}\rho),p)=1$ and thus
$v_L((\gamma-1)\mbox{Tr}_{N/L}\rho)=v_L(\mbox{Tr}_{N/L}\rho)+b$.
Notice that $\mbox{Tr}_{N/L}\rho_* \equiv
(\gamma-1)\mbox{Tr}_{N/L}\rho\bmod p(\gamma-1)\mbox{Tr}_{N/L}\rho$.
Thus $v_L(\mbox{Tr}_{N/L}\rho_*)= v_L(\mbox{Tr}_{N/L}\rho)+b$ as well.
\qed\end{pf}

We next establish that $b_*$ is independent of our choice of group
generators: that a change from $\langle\gamma,\sigma\rangle$ to
$\langle\gamma^a\sigma^b,\gamma^c\sigma^d\rangle$ does not effect
$b_*$, and so we have the identity $v_N((\Theta'-1)\rho_0)-
v_N(\rho_0)=b_*$ where
$\Theta'=(\gamma^a\sigma^b)(\gamma^c\sigma^d)^{[-(\omega_{\gamma^a\sigma^b,\gamma^c\sigma^d})]}$. Let
$\omega'=(c\omega_{\gamma,\sigma}+d)/(ad-bc)$.  Using the fact that
$G^{\mathbb{F}}$ is a vector space over $\mathbb{F}_q$, we have
$(\Theta')^{[\omega']}=
((\gamma^a\sigma^b)^{[c\omega_{\gamma,\sigma}+d]}(\gamma^c\sigma^d)^{[
-(a\omega_{\gamma,\sigma}+b)]})^{[1/(ad-bc)]}=\Theta$ in
$G^{\mathbb{F}}$. And so $(\Theta')^{[\omega']}\in\Theta
(1+\mathcal{A}^p)$.  The following lemma allows us to ignore terms in
$\mathcal{A}^p$, so
$v_N(((\Theta')^{[\omega']}-1)\rho_0)-v_N(\rho_0)=b_*$. The desired
identity then follows since $((\Theta')^{[\omega']}-1)\rho_0\equiv
\omega'(\Theta'-1)\rho_0\bmod (\Theta'-1)\rho_0\euP_N$.

\begin{lem} 
Given $\mu\in\mathcal{A}^p$ or $\mu\in (\sigma-1)^p\subseteq
\euO_L[\sigma]$, 
then for all $\rho\in N$,
$$v_N(\mu\rho)>v_N(\rho)+pb.$$  In
particular, when $v_{N}(\rho)\equiv b \bmod p$ and
$\kappa_i\in \euO_L$, we have
$$v_N((\sigma^{[\kappa_1+\kappa_2]}-\sigma^{[\kappa_1]}\sigma^{[\kappa_2]})\rho)>v_N(\rho)+pb,$$
and so if
$v_{N}((\gamma\sigma^{[\kappa_1+\kappa_2]}-1)\rho) =v_N(\rho)+pb$, and
$\min\{v_{N}((\gamma\sigma^{[\kappa_1]}-1)\rho),
v_{N}({\kappa_2}(\sigma-1)\rho)\}<v_N(\rho)+pb$, then
$v_{N}((\gamma\sigma^{[\kappa_1]}-1)\rho)=
v_{N}({\kappa_2}(\sigma-1)\rho)$.
\end{lem}

\begin{pf}
Since $v_N((\sigma-1)\rho)\geq v_N(\rho)+b$ and
$v_N((\gamma-1)\rho)\geq v_N(\rho)+b$ with strict inequality when
$v_N(\rho)\equiv 0\bmod p$, we have $v_N(\mu\rho)>v_N(\rho)+pb$ for
all $\rho\in N$.  To prove the rest of the lemma, we need
$\sigma^{[\kappa_1]}\cdot\sigma^{[\kappa_2]}\equiv
\sigma^{[\kappa_1+\kappa_2]}\bmod (\sigma-1)^p$ in $\euO_L[G]$.
So observe that $(1+X)^{Y}\cdot (1+X)^{Z}=(1+X)^{Y+Z}$ in the polynomial
ring $\bQ[X,Y,Z]/(X^p)$. Therefore $(1+X)^{[Y]}\cdot (1+X)^{[Z]}=
(1+X)^{[Y+Z]}$ in $\bZ_{(p)}[X,Y,Z]/(X^p)$. Now set $X=\sigma-1$ to obtain the second statement.  As a
result, if $v_{N}((\gamma\sigma^{[\kappa_1+\kappa_2]}-1)\rho)
-v_N(\rho)=pb$, we have
$(\gamma\sigma^{[\kappa_1]}\sigma^{[\kappa_2]}-1) \rho\equiv 0\bmod
\rho\euP_N^{pb}$ and $(\gamma\sigma^{[\kappa_1]}-1) \rho\equiv
-\gamma\sigma^{[\kappa_1]}\cdot (\sigma^{[\kappa_2]}-1) \rho\bmod
\rho\euP_N^{pb}$.  Since $\gamma\sigma^{[\kappa_1]}$ is a unit and
$(\sigma^{[\kappa_2]}-1)=
\kappa_2(\sigma-1)+\sum_{i=2}^{p-1}\binom{\kappa_2}{i}(\sigma-1)^i$,
the last statement follows.\qed\end{pf}

Our final technical lemma establishes that the value of $b_*$ is
independent of our choice of $\rho_0$.

\begin{lem} Given $\rho\in N$ with 
$v_{N}(\rho)\equiv b \bmod p^2$ and $\kappa\in \euO_L$, let
$\mathcal{B}:=v_{N}((\gamma\sigma^{[\kappa]}-1)\rho)-v_N(\rho)$.  Then for all
$\rho'\in N$, and $\mu\in\mathcal{A}^p$
\[v_{N}((\gamma\sigma^{[\kappa]}(1+\mu)-1)\rho')-v_N(\rho')\geq \mathcal{B}.\]
Moreover, we
have equality in the following cases:
\begin{enumerate}
\item[{\rm (i)}] $\mathcal{B}=pb$,
$v_{N}(\rho')\equiv b \bmod p$, but $v_{N}(\rho')\not\equiv (1-p)b
\bmod p^2$,
\item[{\rm (ii)}] $\mathcal{B}<pb$
and $v_{N}(\rho')\equiv b \bmod p^2$.
\item[{\rm (iii)}] $\mathcal{B}\equiv b\bmod p$
and $v_{N}(\rho')\not\equiv 0 \bmod p$.
\end{enumerate}
\end{lem}

\begin{pf}
Write $(\gamma\sigma^{[\kappa]}(1+\mu)-1)\rho'=A+B$
where $A=\gamma\sigma^{[\kappa]}\mu\rho'$ and
$B=(\gamma\sigma^{[\kappa]}-1)\rho'$. By Lemma 2,
$v_N(A)>v_N(\rho')+pb$. And so by Lemma 1,
$v_N(A)>v_N(\rho')+\mathcal{B}$.  We are left to 
prove $v_{N}(B)\geq v_N(\rho')+\mathcal{B}$, with equality in cases (i)--(iii).

We express $\rho'$ in terms of $\rho$.  Notice that since
$\{v_N((\sigma-1)^i\rho): i=0, \ldots, p-1\}$ is a complete set of
residues modulo $p$ and $N/L$ is fully ramified, there are $a_i\in L$
such that $\rho'=\sum_{i=0}^{p-1}a_i(\sigma-1)^i\rho$.  Choose $i_0$
such that $v_N(\rho')=v_N(a_{i_0})+i_0b +v_N(\rho)\equiv (i_0+1)b\bmod
p$. So $v_N(a_{i_0}(\sigma-1)^{i_0}\rho)=v_N(\rho')$.  Note that for
$i\neq i_0$, $v_N(a_i(\sigma-1)^i\rho)>v_N(\rho')$ and so $v_N(a_i)+ib
>v_N(a_{i_0})+i_0b$.  For each $i$, let $(\gamma\sigma^{[\kappa]}-1)
\cdot a_i(\sigma-1)^i\rho = A_i +B_i$ where $A_i=[(\gamma-1)a_i]\cdot
\gamma\sigma^{[\kappa]}(\sigma-1)^i\rho$ and $B_i=
a_i(\sigma-1)^i\cdot (\gamma\sigma^{[\kappa]}-1)\rho$. This means that
$B=\sum_{i=0}^{p-1}(A_i+B_i)$.
Our goal is to prove that $v_N(\sum_{i=0}^{p-1}(A_i+B_i))\geq
v_N(\rho')+\mathcal{B}$.

Begin with the $A_i$. Notice that since $\gamma\sigma^{[\kappa]}$ is a
unit, $v_N(A_i)= v_N((\gamma-1)a_i)+ v_N((\sigma-1)^i\rho)$, where
$v_N((\gamma-1)a_i)\geq v_N(a_i)+pb$.  So for $i\neq i_0$, we have
strict inequality, $v_N(A_i)>v_N(\rho')+pb\geq
v_N(\rho')+\mathcal{B}$.  For $i=i_0$, we have $v_N(A_{i_0})\geq
v_N(\rho')+pb\geq v_N(\rho')+\mathcal{B}$ with strict inequality when
$\mathcal{B}<pb$.

Consider the $B_i$. Note that since $v_N(\rho')=v_N(a_{i_0})+i_0b
+v_N(\rho)$, we have $v_N(B_{i_0})\geq v_N(a_{i_0})+i_0b+
v_N((\gamma\sigma^{[\kappa]}-1)\rho) =v_N(\rho')+\mathcal{B}$.  For
$i\neq i_0$ we have $v_N(a_i)+ib>v_N(a_{i_0})+i_0b$, and so we have
strict inequality $v_N(B_i)> v_N(\rho')+\mathcal{B}$.

When do we have equality in the statement of our lemma? Case (i) is
clear and follows immediately from Lemma 1. In cases (ii) and (iii) we
have $\mathcal{B}<pb$, and so equality occurs precisely when
$v_N(B_{i_0})= v_N(\rho')+\mathcal{B}$, which occurs if and only if
$v_N((\sigma-1)^j(\gamma\sigma^{[\kappa]}-1)\rho)\not\equiv 0\bmod p$
for each $0\leq j\leq i_0-1$. There are two extreme cases where this
condition is easy to check. when $i_0=0$, the condition is empty. This
is case (ii). When $\mathcal{B}\equiv b\bmod p$, we have
$v_N((\gamma\sigma^{[\kappa]}-1)\rho)\equiv 2b\bmod p$ and so
$v_N((\sigma-1)^j(\gamma\sigma^{[\kappa]}-1)\rho)\not\equiv 0\bmod p$
for $0\leq j\leq p-3$. The condition holds then if $i_0\leq p-2$,
which is equivalent to $v_{N}(\rho')\not\equiv 0 \bmod p$. This is
case (iii). \qed\end{pf}

Based upon these technical results, the integer $b_*$ satisfies
$b<b_*\leq pb$ and is canonical (independent of our choice of group
generators and element $\rho\in N$ with $v_N(\rho)\equiv b\bmod p^2$).
This is collected in the following theorem where we prove that it is also
the second refined break, as defined in \S1.1.

\begin{thm}
Let $K$ be a finite extension of the field $\bQ_p$ of $p$-adic numbers
with absolute ramification index $e_K$ and inertia degree $f$.
Let $N/K$ be a fully ramified, bicyclic extension with one
ramification break at $b$, and let
$G=\mbox{Gal}(N/K)=\langle\gamma,\sigma\rangle$. Pick any $\rho\in N$
with $v_N(\rho)=b$. Define $\omega$ to be the unique $p^f-1$ root of
unity such that $v_N((\gamma-1)\rho+\omega(\sigma-1)\rho)>2b$, and let
$\Theta=\gamma\sigma^{[\omega]}\in\euO_T[G]$.  Then the refined
ramification filtration has two breaks $b$ and
$b_*=v_N((\Theta-1)\rho)-v_N(\rho)$, so that
\[G^{\bF}=\langle\Theta,\sigma\rangle=G^{\bF,
\rho}_b\supsetneq G^{\bF, \rho}_{b+1}= \langle\Theta\rangle=G^{\bF,
\rho}_{b_*}\supsetneq G^{\bF, \rho}_{b_*+1}=\{e\}. \]
Moreover $b_*$ satisfies $b<b_*\leq pb$ and
is independent of our choices..
\end{thm}
\begin{pf} 
By Lemma 1, $b<b_*\leq pb$.  Let $\bar{\Theta}$ denote the image of
$\Theta$ in $\mathcal{G}$ as defined in \S1.1. By Lemma 3 cases (i) and (ii)
(with $\rho'=\rho$), we have
$i_{\rho}(\bar{\Theta})=b+b_*$.
\qed\end{pf}

\subsection{The value of the second refined break}

The determination of all possible values of $b_*$ will require a
detour through (and detailed analysis of) Kummer bicyclic extensions
with one break at $b$. We therefore begin by summarizing the results
of this detour in the following theorem, which is a consequence of
Proposition 10 and Corollary 11.  Its proof appears in \S2.2.4.

\begin{thm} Let $U:=pb-\max\{(p^2-1)b-p^2e_K,0\}$.
Assuming the conditions of Theorem 4,
$b<b_*\leq U$, and
if $b_*< U$ then $b_*\equiv b\bmod p$ but
$b_*\not\equiv (1+p)b\bmod p^2$.
Moreover, any integer
that satisfies these conditions is the second refined break of
a bicyclic extension with one break at $b$.
\end{thm}
\begin{cor}
$$pb_*-b<p^2e_K$$
\end{cor}
\begin{pf} From Theorem 5,
$b_*\leq pb-\max\{(p^2-1)b-p^2e_K,0\}$, which leads to two cases
depending upon whether or not $\max\{(p^2-1)b-p^2e_K,0\}=0$. Suppose
$\max\{(p^2-1)b-p^2e_K,0\}=0$. 
Thus $(p^2-1)b<p^2e_K$ (recall $\gcd(p,b)=1$)
and $b_*\leq pb$.  Then $pb_*-b\leq (p^2-1)b<p^2e_K$.  Suppose
$\max\{(p^2-1)b-p^2e_K,0\}=(p^2-1)b-p^2e_K$.  Thus $p^2e_K<(p^2-1)b$
and $b_*\leq p^2e_K-(p^2-p-1)b$. So $pb_*-b\leq
p^3e_K-(p^2-1)(p-1)b<p^3e_K-(p-1)p^2e_K=p^2e_K$.  \qed\end{pf}

\subsubsection{A brief history}
The chronology of this research may be of interest.
We began our investigations by
looking at Kummer extensions, as we tried to generalize the results of
\cite{elder:onebreak} from $p=2$ to $p>2$. 
In the course of these investigations, truncated
exponentiation appeared first within the group ring
$\euO_T[G]$, as we worked to prove Lemma 9. It is this appearance of
truncated exponentiation that led us to the investigations in
\cite{elder:newbreaks}, and to Lemmas 1, 2 and 3 and Theorem 4.  Only
later as we worked to determine the precise value of $b_*$, did
truncated exponentiation emerge among the generators of the bicyclic
Kummer extension. This work is captured in Proposition 10 below.  Our
presentation here reverses that chronology somewhat, as we start in
\S2.2.2 by assuming truncated exponentiation among the generators of
our extension.

\subsubsection{Bicyclic Kummer extensions with one break}

Let $\zeta$ denote a nontrivial $p$th root of unity, and assume that
$\zeta\in K$.  Given any integer $b$ such that $0< b<pe_K/(p-1)$ with
$\gcd(b,p)=1$, choose $\beta\in K$ such that
$v_K(\beta)=pe_K/(p-1)-b$.  Choose a $p^f-1$ root of unity
$\omega$ such that $\omega^{p-1}\neq 1$, and set
\[x^p=1+\beta, \quad
y^p=(1+\beta)^{[\omega^p]}.\]
For either $t=0$ or $0<t<b$ with $\gcd(t,p)=1$, choose
$\tau\in K$ such that $v_K(\tau)=pe_K/(p-1)-t$.  Set
\[z^p=1+\tau.\] 

Then $N_z:=K(x,yz)$, a subfield of $K(x, y, z)$, is a fully ramified,
bicyclic extension with one break in its ramification filtration, at
$b$.  Moreover, {\em any} fully ramified, bicyclic extension with one
break can be represented in this way. In particular, there are $\tau$
with $t=0$ such that $1+\tau$ is a $p$th power. In this case, we
have $N_z=N_1:=K(x,y)$.

Choose $\sigma, \gamma\in
G=\mbox{Gal}(N_z/K)$ with 
\[\begin{array}{ll}
\sigma x=x, & \sigma yz=\zeta yz, \\
\gamma x=\zeta x,& \gamma yz= yz.
\end{array}\]
And let $L=K(x)$.

Why have we chosen to express the generators in this way?  Our first
choice, to represent $x^p$ as $1+\beta$, is natural: $p$-adic defects
of units are related to ramification numbers \cite{wyman}.  Our second
choice, to represent $yz$ as a product, means that $N_z$ can be seen
as a `twist' of $N_1=K(x,y)$. See \S2.2.3.  Our final choice, to
relate $y^p$ to $x^p$ by truncated exponentiation, is justified simply
by the fact that it makes the nice statement in Proposition 10
possible.

We are interested in the refined ramification filtration, and so we
require now an element $\rho_0$ of $N_z$ with valuation $b$. Observe
that since $N_z/L$ is a cyclic Kummer extension with break number $b$,
$N_z=L(Y_z)$ for some $Y_z$ with $Y_z^p=1+\beta_z\in L$ and
$v_L(\beta_z)=p^2e_K/(p-1)-b$.  Clearly then
$\rho_0=(\zeta-1)/(Y_z-1)$ will do.  Observe furthermore
$L(Y_z)=L(yz)$.  To describe the Galois action (and in particular the
$\gamma$-action) on $\rho_0$ and thus on $Y_z$ we ask that $yz/Y_z$ be
an explicitly described element in $L$. This is accomplished in the
following two lemmas.

\begin{lem}
There is a $\beta_L\in L$ with $v_L(\beta_L)=p^2e_K/(p-1)-b$ such that
\[(1+\beta)^{[\omega^p]}= \left(x^{[\omega]}\right)^p\cdot (1+\beta_L).\]
\end{lem}
\begin{pf}
The norm, from $L$ to $K$, of $x-1$ is $(-1)^{p-1}\beta$.
So $v_L(x-1)=pe_K/(p-1)-b$ and thus $v_L(p(x-1))=p^2e_K/(p-1)-b$.  Now
$\left(x^{[\omega]}\right)^p\equiv
\sum_{i=0}^{p-1}\binom{\omega}{i}^p(x-1)^{pi} +\omega p(x-1)\bmod
p(x-1)^2$. Note that $\binom{\omega}{i}^p=\binom{\omega^p}{i}$ for
$i=0, 1$, and $\binom{\omega}{i}^p(x-1)^i\equiv
\binom{\omega^p}{i}(x-1)^i\bmod p(x-1)^2$ for $i\geq 2$.  Furthermore
since $1+\beta=(1+(x-1))^p$,
$(x-1)^p=\beta-\sum_{i=1}^{p-1}\binom{p}{i}(x-1)^i\equiv
\beta-p(x-1)\bmod p(x-1)^2$. So $(x-1)^{pi}\equiv \beta^i\bmod
p(x-1)^2$ for $i>1$. Therefore 
$\left(x^{[\omega]}\right)^p\equiv
\sum_{i=0}^{p-1}\binom{\omega^p}{i}\beta^i\cdot\left(1
+(\omega-\omega^p) p(x-1)\right)\bmod p(x-1)^2$.  Since
$\omega\not\in\bZ_p$, $\omega-\omega^p$ is a unit. The result follows.
\qed\end{pf}

\begin{lem} 
There are elements $\delta', \tau_L\in L$ with
$v_L(\delta')=pe_K/(p-1)-t$ and $v_L(\tau_L)=p^2e_K/(p-1)-t$ such that
$1+\tau=(1+\delta')^p(1+\tau_L)$.
\end{lem}
\begin{pf} 
If $t=0$ then $K(z)/K$ is unramified. Thus $L(z)/L$ is unramified and
the result is clear. If $t\neq 0$ then $K(z)/K$ is ramified with
ramification number $t$. Thus $K(x,z)$ is a fully ramified $C_p^2$
extension with two lower ramification numbers, $b_1=t,
b_2=t+p(b-t)$. Since $L(z)/L$ is a Kummer, ramified $C_p$-extension
with ramification number $t$, we find that $L(z)=L(Z)$ where
$Z^p=1+\tau_L$ for some $\tau_L\in L$ with
$v_L(\tau_L)=p^2e_K/(p-1)-t$ \cite{wyman}.  Moreover, $Z$ may be
chosen so that $z/Z\in L$. In that case, $z/Z=1+\delta'$ for some
$\delta'\in L$ with $v_L(\delta')=pe_K/(p-1)-t$.
\qed\end{pf}

Now using the $\delta'$ of Lemma 8, define $r_z\in L$ by 
\begin{equation}r_z= x^{[\omega]}(1+\delta)\mbox{ where }
\delta= \begin{cases}\delta' \mbox{
for $t>b/p$,}\\ 0\mbox{ for $t<b/p$.}
\end{cases}
\end{equation}

Choose $Y_z=yz/r_z\in N_z$, so $r_z$ is the `ratio'
$yz/Y_z\in L$ and $\sigma Y_z=\zeta Y_z$.  Using Lemma 7, 
$Y_z^p=1+\beta_z$ where
\[1+\beta_z=
\begin{cases} (1+\beta_L)(1+\tau_L) \mbox{ for $t>b/p$},\\ 
(1+\beta_L)(1+\tau) \mbox{ for $t<b/p$}.
\end{cases}\]
As a result, $v_{N_z}(Y_z-1)= v_L(\beta_z)= p^2e_K/(p-1)-b$ and
\begin{equation}
\rho_0=\frac{\zeta-1}{Y_z-1}
\end{equation} 
satisfies $v_{N_z}(\rho_0)=b$.\

We now recall an earlier observation: Since
$v_{N_z}((\gamma-1)\rho_0)=v_{N_z}((\sigma-1)\rho_0)=2b$, there is an
element $a\in \euO_T$ such that $(\gamma-1)\rho_0\equiv
a(\sigma-1)\rho_0\bmod \euP_{N_z}^{2b+1}$, which can be rewritten as
$\gamma\rho_0\equiv\sigma^{[a]}\rho_0\bmod \euP_{N_z}^{2b+1} $, and
also as $(\gamma\sigma^{[-a]}-1)\rho_0\bmod \euP_{N_z}^{2b+1}$. We are
interested in determining $a$ along with the precise valuation,
$v_{N_z}((\gamma\sigma^{[-a]}-1)\rho_0)$. Recall the generic bounds
given in Lemma 1.

\begin{lem} Using the notation of this section, $\gamma\sigma^{[-\Omega_z]}\rho_0\equiv \rho_0 \bmod
\rho_0^{1+p}$ where
\[ \Omega_z:=\frac{(\gamma-1)Y_z}{(\sigma-1)Y_z}\in \euO_L^*.
\]
\end{lem}
\begin{pf}
Using the fact that $\sigma Y_z=\zeta Y_z$, we find that
\[\sigma\rho_0=\frac{\rho_0}{1+Y_z\rho_0}\equiv
\frac{\rho_0}{1+\rho_0} \bmod (\zeta-1)\rho_0.\] So we can establish
by induction that
\[(\sigma-1)^t\rho_0\equiv(-1)^t t!\prod_{i=0}^t\frac{\rho_0}{1+i\rho_0} \bmod \rho_0(\zeta-1)\mbox{ for } 0\leq t\leq p-1.\]
Now define $[X]_n=X(X-1)\cdots (X-(n-1))\in \bZ[X]$ so that
$\binom{X}{n}\cdot n!=[X]_n$ and establish the following power series
identity for $\Omega\in \euO_L$ by induction 
\[\sum_{s=0}^t(-1)^s[\Omega]_s\prod_{i=1}^s\frac{X}{1+iX}=\frac{1}{1+\Omega X}\left (1+(-1)^t[\Omega]_{t+1}\prod_{i=0}^t\frac{X}{1+iX}\right )\in \euO_L[[X]].\]
As a result,
\[\sigma^{[\Omega]}\rho_0\equiv \frac{\rho_0}{1+\Omega\rho_0}\left
(1+(\Omega^p-\Omega)\frac{\rho_0^p}{1-\rho_0^{p-1}}\right )\bmod
(\zeta-1)\rho_0\]
and thus
\[\sigma^{[\Omega]}\rho_0\equiv \frac{\rho_0}{1+\Omega\rho_0}\bmod \rho_0^{1+p}.\]
Now observe that since $v_{N_z}((\gamma-1)Y_z)= v_{N_z}(\zeta-1)=
v_{N_z}((\sigma-1)Y_z))$,
\[\gamma\rho_0=\frac{\rho_0}{1+\frac{(\gamma-1)Y_z}{\zeta-1}\rho_0}
\equiv \frac{\rho_0}{1+\Omega_z\rho_0} \bmod (\zeta-1)\rho_0\] where
$\Omega_z$ is as above.
Putting these together yields $\sigma^{[\Omega_z]}\rho_0\equiv
\gamma\rho_0\bmod \rho_0^{1+p}$. By Lemma 2,
$\sigma^{[-\Omega_z]}\sigma^{[\Omega_z]} \rho_0\equiv \rho_0\bmod
\rho_0^{1+p}$. Thus the desired statement holds.\qed\end{pf}
\begin{prop} 
$\Omega_z\equiv -\omega\bmod \euP_L$. Thus $b_*=v_{N_z}((\Theta-1)\rho_0)-v_{N_z}(\rho_0)$ where $\Theta=\gamma\sigma^{[\omega]}$.
Let $\eta_z:=\Omega_z+\omega\in \euP_L$. Then for $b_*<pb$,
\[v_L(\eta_z)=\frac{b_*-b}{p}.\]
In general,
$b_*=pb-\max\{(p^2-1)b-p^2e_K, pt-b, 0\}$.
\end{prop}

\begin{pf}
Recall the unit $r_z$. Using its definition in (1), we find that
$(\gamma-1) r_z=((\zeta x)^{[\omega]}- x^{[\omega]})(1+\delta) +
(\zeta x)^{[\omega]}((\gamma-1)\delta)$.  Our first observation is
that since $v_L((\gamma-1)\delta)\geq pe_K/(p-1)-t+b> v_L(\zeta-1)$,
we have $(\gamma-1) r_z\equiv 0\bmod (\zeta-1)$.  So using
$Y_z=yz/r_z$, we can decompose $\Omega_z$ as a product:
$\Omega_z=-A\cdot B$ with $A:=(\gamma r_z)^{-1}\equiv r_z^{-1} \bmod
(\zeta-1)$ and $B:= (\gamma-1)r_z/(\zeta-1)\in \euO_L$.

To describe $B$ further, we examine the term $C:=(\zeta x)^{[\omega]}-
x^{[\omega]}$ modulo $(\zeta-1)^2$. For $1\leq i\leq p-1$,
$(\zeta x-1)^i=((\zeta-1)x+(x-1))^i\equiv
i(\zeta-1)x(x-1))^{i-1}+(x-1)^i\bmod (\zeta-1)^2$. So $C\equiv
(\zeta-1)\sum_{i=1}^{p-1}\binom{\omega}{i}ix(x-1)^{i-1}\bmod
(\zeta-1)^2$. Observe that
$\binom{\omega}{i}i=\omega\binom{\omega-1}{i-1}$. So $C\equiv
(\zeta-1)\cdot \omega x\left
[x^{[\omega-1]}-\binom{\omega-1}{p-1}(x-1)^{p-1}\right ]\bmod
(\zeta-1)^2$. Now replace $A$, $B$ and $C$, in the expression for
$\Omega_z$, and find
\begin{multline*}
\Omega_z\equiv-\frac{\omega x \left
(x^{[\omega-1]}-\binom{\omega-1}{{p-1}}(x-1)^{p-1}\right )}
{x^{[\omega]}} -\frac{(\gamma-1)\delta}{\zeta-1}\frac{1}{1+\delta}
\bmod(\zeta-1)\\ \equiv-\omega +\omega\binom{\omega-1}{{p-1}}(x-1)^{p-1}
-\frac{(\gamma-1)\delta}{\zeta-1} \bmod\left
(\zeta-1,(x-1)^{p}, \delta\frac{(\gamma-1)\delta}{\zeta-1} \right ),
\end{multline*}
which proves the first assertion and establishes a congruence relation
$\eta_z\equiv D-E$ with $D:= \omega\binom{\omega-1}{{p-1}}(x-1)^{p-1}$
and $E:= (\gamma-1)\delta/(\zeta-1)$.  There are two cases to
consider: $t<b/p$ and $t>b/p$. If $t<b/p$, $\delta=0$ and so $E=0$.
Since $v_L(D)<v_L((x-1)^p$, we find $v_L(\eta_z)=v_L(D)=pe_K-(p-1)b$,
when $v_L(\eta_z)<v_L(\zeta-1)$.  On the other hand, if $t>b/p$ then
we have $\gcd(t,p)=1$ and since $v_L(D)=pe_K-(p-1)b\equiv b\bmod p$
while $v_L(E)=b-t$, $v_L(D)\not\equiv v_L(E)\bmod p$. Thus
$v_L(D-E)=\min\{v_L(D),v_L(E)\}$. Since $\min\{v_L(D),v_L(E)\} <
\min\{v_L(D(x-1)),v_L(\delta E)\}$, we have
$v_L(\eta_z)=\min\{v_L(D),v_L(E)\}=b-\max\{p(b-e_K),t\}$, whenever
$v_L(\eta_z)<v_L(\zeta-1)$.

From Lemma 1,
$b_*:=v_{N_z}((\gamma\sigma^{[\omega]}-1)\rho_0)-v_{N_z}(\rho_0)\leq pb$.
Using Lemma 2 with $\kappa_1=\omega$ and $\kappa_2=-\eta_z$, we find that when
$b_*<pb$, we have
$v_{N_z}((\gamma\sigma^{[\omega]}-1)\rho_0)=v_{N_z}(\eta_z(\sigma-1)\rho_0)$.
So $b_*=pv_L(\eta_z)+b$. But then $v_L(\eta_z)<b-b/p<v_L(\zeta-1)$. So
substituting our formulas for $v_L(\eta_z)$ into $b_*=pv_L(\eta_z)+b$,
we find that $b_*=pb-(p^2-1)b+p^2e_K$ for $t<b/p$ and
$b_*=pb-\max\{(p^2-1)b-p^2e_K),pt-b\}$ for $t>b/p$.  Since $b_*\leq
pb$, these both agree with $b_*=pb-\max\{(p^2-1)b-p^2e_K),pt-b,0\}$.
\qed\end{pf}

\begin{cor}
Let $U:=pb-\max\{(p^2-1)b-p^2e_K,0\}$. 
Any integer $n$ satisfying $b<n\leq U$,
and if $n< U$ then $n\equiv b\bmod p$ but
$n\not\equiv (1+p)b\bmod p^2$, is the second refined break for
a bicyclic Kummer extension with one break at $b$.
\end{cor}

\subsubsection{Strong twists alter ramification breaks}

Let $\bar{G}=\mbox{Gal}(\bar{K}/K)$ denote the absolute Galois group.
We will call the fixed field of the kernel of a representation of
$\bar{G}$, the fixed field of the representation.  Let $\chi_x$,
$\chi_{xy}$, $\chi_z$ be $1$-dimensional characters with fixed fields
$K(x)$, $K(xy)$ and $K(z)$ respectively.  Let $V$ denote the
$2$-dimensional representation of $\bar{G}$ with character
$\chi_y+\chi_{xy}$ and fixed field $N_1=K(x,y)$. Then $N_z=K(x,yz)$ is
the fixed field of the twisted representation $V\otimes \chi_z$.  The
`strength' of the twist by $\chi_z$ is parametrized by $t$, the
ramification break of $K(z)/K$.

Consider the following diagram with the formula for $b_*$ displayed as
a function of $(b,t)$ in each of three relevant regions that lie
below the diagonal line $t=b$.  The boundaries of these regions are:
the line $t=b$; the segment $\ell_1$, which is on the line $t=b/p$;
the segment $\ell_2$, which is on $t=p(b-e_K)$; and the segment
$\ell_3$, which is on $b=p^2e_K/(p^2-1)$.


\hspace*{.75in}\input{perm-graph}

Now view $N_z$ is a twist of $N_1$ and observe that `strong' twists
change ramification filtrations, while `weak' twists preserve them: If
the twist is `weak' and thus $t$ is relatively small ($t<b/p$ or
$t<p(b-e_K)$), the formula for $b_*$ in $N_z$ is the same as in
$N_1$. Otherwise the formulas for $b_*$ are different (although if
$t<b$, $N_z/K$ still has only one ramification break).  If we strengthen
our twist further and choose $t>b$, then $N_z/K$ will have two ramification
breaks.

Why is this so? Why are the values of the second refined breaks in
$N_z$ and $N_1$ equal when $t<b/p$ or $t<p(b-e_K)$?  Observe that the
formula for $b_*$ results from the expression for $v_L(\eta_z)$
determined in Proposition 10. Note furthermore that the proof of
Proposition 10 describes $\eta_z$ completely in terms of $r_z$.  So
our question becomes: Why do $r_z,r_1\in L$ ``agree'' under $t<b/p$ or
$t<p(b-e_K)$? When $t<b/p$, because they are equal. Recall (1).  So
where it matters, the twist has no effect!  Now consider $t<p(b-e_K)$
with $t>b/p$. Motivated by our answer for $t<b/p$, observe that
$t<p(b-e_K)$ is equivalent to $v_L(\beta)<v_L(\delta')$, where
$\delta'$ was defined in Lemma 8. Returning to (1), we conclude that
they ``agree'' because they are equivalent, $r_z\equiv r_1\bmod
\beta\euP_L$.

\subsubsection{Bicyclic non-Kummer extensions with one break}
\begin{pf*}{Proof (Theorem 5)}
Recall that $N/K$ is a fully ramified, bicyclic extension with one
ramification break at $b$. If $\zeta\in K$ and so the $p$th roots of
unity are present, the result is contained in Proposition 10 and
Corollary 11. To apply these results when $\zeta\not\in
K$, we consider the related Kummer extension $N(\zeta)/K(\zeta)$ with
$d=[K(\zeta):K]$.  By abuse of notation use $\sigma, \gamma$ to
represent automorphisms in $\mbox{Gal}(N(\zeta)/K)$, so that
$\langle\sigma, \gamma\rangle=
\mbox{Gal}(N(\zeta)/K(\zeta))=\mbox{Gal}(N/K)$. Pick any $\rho_0\in N$
with $v_N(\rho_0)=b$. Then $v_{N(\zeta)}(\rho_0)=db$. Using the
Herbrand function \cite[IV \S3]{serre:local}, the ramification break
of $N(\zeta)/K(\zeta)$ is $db$.  Recall from the beginning of \S2,
that $\omega$ is defined to be the unique $p^f-1$ root of unity such
that $(\gamma-1)\rho_0\equiv
-\omega(\sigma-1)\rho_0\bmod\euP_N^{2b+1}$.
Since
$(\gamma\sigma^{[\omega]}-1)\rho_0\equiv 0\bmod\rho_0^2\pi_N$ in $N$,
$(\gamma\sigma^{[\omega]}-1)\rho_0\equiv 0\bmod\rho_0^2\pi_{N(\zeta)}$
in $N(\zeta)$. Therefore the $\omega$ defined here is the
same as the $\omega$ defined in \S2.2.2 for $N(\zeta)/K(\zeta)$.  And
$v_{N(\zeta)}((\Theta-1)\rho_0)=db +db_*$, where $db_*$ is determined
by Proposition 10 with $b$ replaced by $db$ and
$e_{K(\zeta)}=de_K$. The result follows now after the integer $d$ is
removed everywhere.  \qed\end{pf*}

\section{Galois Module Structure in Bicyclic Extensions}

We are interested in the relevance of the second refined break $b_*$
for Galois module structure.  Let $N/K$ be a fully ramified, bicyclic
extension with one break $b$ in its ramification filtration and assume
the notation of \S2.

In Theorem 12 of \S3.1, we determine just enough of the $\bF_q[G]$-structure of
$\euP_N^r/p\euP_N^r$  to prove that, if $b_*<(p-1+1/p
)b$, this structure depends upon $b_*$.  As a result, the
$\euO_T[G]$-structure of ideals also depends upon $b_*$.

Next, because it is easily done, we assume in \S3.2 that we have
maximal refined ramification $b_*=pb$, and in Theorem 18 explicitly
describe, in a transparent way, the $\euO_T[G]$-structure of each
ideal $\euP_N^r$.

Based upon \cite{elder:onebreak}, we conjecture that the our result
concerning the relevance of $b_*$ is sharp -- namely, that the
$\euO_T[G]$-structure of each ideal $\euP_N^r$ under $(p-1+1/p
)b<b_*<pb$, which we call {\em near maximal refined ramification}, is
independent of $b_*$ and in fact agrees with the structure given in
Theorem 18.

\subsection{On modular Galois module structure}

Identify $\Theta\in\euO_T[G]$ with its image in $\bF_q[G]$, and
observe that $(\Theta)^p=1$ in $\bF_q[G]$. There are exactly $p$
indecomposable modules over $\mathbb{F}_q[\Theta]$, namely
$L(i)=\bF_q[x]/(x-1)^i$ for $1\leq i\leq p$, where $\Theta$ acts via
multiplication by $x$. This means that $\euP_N^r/p\euP_N^r$ is
uniquely expressible as
$$\euP_N^r/p\euP_N^r \cong\bigoplus_{i=1}^p L(i)^{a_i}$$ for some
integers $a_i\geq 0$. Here we determine $a_p$, and in particular find

\begin{thm}
\begin{multline*}
a_p=\dim_{\bF_q}\left ((\Theta-1)^{p-1}\euP_N^r/p\euP_N^r\right )
=\\pe_K+\left\lceil\frac{r}{p}\right\rceil
-\left\lceil\frac{r-b}{p}\right\rceil
-b-\begin{cases}
\frac{b_*-b}{p} \mbox{ for
}b_*<(p-1+1/p)b,  \vspace*{3mm}\\
b+\left\lceil\frac{r-pb}{p^2}\right\rceil-\left\lceil\frac{r+(p-1)b}{p^2}\right\rceil
\mbox{ otherwise}.
\end{cases} 
\end{multline*}
\end{thm}

This result for $b_*=pb$ follows from Theorem 18.
In this section, we verify it for $b_*<pb$, which allows
us to use the fact that $c=(b_*-b)/p$ is an integer.

We begin by establishing an $\euO_T$-basis for $\euP_N^r$, a basis
that will also serve as a $\bF_q$-basis for
$$\mathcal{M}=\euP_N^r/p\euP_N^r.$$ 
Let $\rho_m\in N$ be any element with $v_N(\rho_m)=b+pm$ and
observe that since $b_*\equiv
b\bmod p$ and $\gcd(b,p)=1$, $\{v_N((\Theta-1)^ip\rho_m): i=0, \ldots,
p-1\}$ is a complete set of residues modulo $p$.  As $m$ varies over
$\bZ$, the resulting elements $(\Theta-1)^ip\rho_m$ will lie in
one-to-one correspondence, via valuation $v_N$, with $\bZ$.  Collect
those with $r\leq v_N((\Theta-1)^ip\rho_m) \leq r+p^2e_K-1$.  We have
a $\euO_T$-basis for $\euP_N^r$. So that we can follow the effect of
$\Theta$ upon this basis, we will replace certain $\rho_m$ with
$\rho_m^*$ of equal valuation. This is done in Lemma 14. But first 
we require a technical lemma.

\begin{lem}
For any $\omega \in \euO_T$, we have the congruence in $\euO_T[\sigma]$
$$ (\sigma^{[\omega]})^p -1 \equiv (w-w^p) \sum_{i=1}^{p-1}
\binom{p}{i} (\sigma-1)^i \bmod{p^2\euO_T[\sigma]}. $$
Recall $\Theta=\gamma\sigma^{[\omega]}\in\euO_T[G]$. Then there is a
unit $u(\sigma)\in\euO_T[\sigma]^*$ defined by
\[(\Theta-1)^p + \sum_{i=1}^{p-1}\binom{p}{i}
(\Theta-1)^i =p(\sigma-1)u(\sigma),\]
satisfying
\[u(\sigma) \equiv
(\omega-\omega^p)\sum_{i=1}^{p-1} \left[ \frac{1}{p} \binom{p}{i}
\right] (\sigma-1)^{i-1} \equiv (\omega-\omega^p)\sum_{i=1}^{p-1}
(-1)^{i+1}i^{-1} (\sigma-1)^{i-1}\]
modulo ${p\euO_T[\sigma]}$.
In particular, 
$(\Theta-1)^p
=[u(\sigma)(\sigma-1)-w(\Theta)(\Theta-1)]\cdot p$
where $u(\sigma)$ and
$w(\Theta)=\sum_{s=1}^{p-1}p^{-1}\binom{p}{s}(\Theta-1)^{s-1}$ are both units
in $\euO_T[G]$.
\end{lem}
\begin{pf}
We work initially in the truncated polynomial ring
$\bQ[W,F]/(F^{2p})$. In this ring we have the (finite) binomial
expansion $ (1+F)^W = \sum_{i=0}^{2p-1} \binom{W}{i} F^i=
\sum_{i=0}^{p-1} \binom{W}{i} F^i + \sum_{i=0}^{p-1} \binom{W}{p+i}
F^{p+i}$.  Now, as observed in the proof of \cite[Lem
2.2]{elder:newbreaks}, for $0 \leq i \leq p-1$ we have $ p \binom{W}{p+i} \in \bZ_{(p)}[W]$, $p\binom{W}{p+i} \equiv (W-W^p)\binom{W}{i}
\pmod{p\bZ_{(p)}[W]}$.  Hence there is a polynomial $e(F,W) \in
\bZ_{(p)}[W,F]$ such that
\begin{eqnarray*}
(1+F)^W & = & \sum_{i=0}^{p-1}\binom{W}{i} F^i + \frac{F^p}{p} \left(
\sum_{i=0}^{p-1} (W-W^p) \binom{W}{i} F^i +p\cdot e(F,W) \right) \\ &
= & (1+F)^{[W]} \left( 1 + (W-W^p) \frac{F^p}{p} \right) + F^p e(W,F).
\end{eqnarray*}
Raising both sides to the power $p$, using $(F^p)^2=0$, and
observing that $( (1+F)^W)^p=( (1+F)^p)^W$ by
properties of (infinite) binomial series, we obtain the following
identity in $\bZ_{(p)}[W,F]/(F^{2p})$:  
$$ \bigl( (1+F)^p \bigr)^W = \bigl( (1+F)^{[W]} \bigr)^p \bigl( 1+
(W-W^p)F^p \bigr) + p\bigl( (1+F)^{[W]} \bigr)^{p-1} F^p
e(W,F). $$ 
Consider its image under the
homomorphism from $\bZ_{(p)}[W,F]/(F^{2p})$ to
$R=(\euO_T/p^2\euO_T)[\sigma]$ which takes $W$ to $\omega$
and $F$ to $f=\sigma-1$. This homomorphism is well-defined
because $f^p=-\sum_{i=1}^{p-1} \binom{p}{i} f^i \in pR$, so
that $f^{2p}\in p^2R$.  In $R$ this identity becomes
\[ 1= (\sigma^{[\omega]})^p\bigl (1+ (\omega-\omega^p)f^p \bigr) .\]
Note that $(1+(\omega-\omega^p)f^p)(1-(\omega-\omega^p)f^p)=1$ in $R$,
and so $(\sigma^{[\omega]})^p= 1- (\omega-\omega^p)f^p $ in $R$.
Moreover expanding $\sigma^p=(1+f)^p$ using the binomial theorem
yields $f^p=-\sum_{i=1}^{p-1} \binom{p}{i}f^i$, and thus
$(\sigma^{[\omega]})^p-1= (\omega-\omega^p)\sum_{i=1}^{p-1}
\binom{p}{i}f^i$ in $R$.

Use the binomial expansion $ (\sigma^{[\omega]})^p=\Theta^p= \bigl(
(\Theta-1) +1\bigr)^p=(\Theta-1)^p+1 + \sum_{i=1}^{p-1}\binom{p}{i}
(\Theta-1)^i$, and we obtain the statements concerning
$u(\sigma)$.  \qed\end{pf}

\begin{lem}
There are $\rho_m, \rho_m^*\in N$ with
$v_N(\rho_m)=v_N(\rho_m^*)=b+pm$ satisfying
\begin{multline*}
(\Theta-1)p\rho_m-
(\Theta-1)^p\rho_{m+c}^*-
\sum_{i=1}^{p-1}
\left[ \frac{1}{p}\binom{p}{i}\right]
(\Theta-1)^ip\rho_{m+c}^*  \\ \vspace*{6mm}
=
\begin{cases}
p\rho_{m+b}&\mbox{for }m\not\equiv -b\bmod p\\
p\rho_{m+pe_K-(p-2)b}&\mbox{for }m\equiv -b\bmod p
\end{cases}
\end{multline*}
\end{lem}
\begin{pf}
For $0\leq k < e_K$ choose $\alpha_k\in L$ with $v_L(\alpha_k)=b+pk$.
Since $u_{\gamma}:=\sum_{i=1}^{p-1}\left
[\frac{1}{p}\binom{p}{i}\right ](\gamma-1)^{i-1}\in\bZ_p[\gamma]^*$,
$v_L(u_{\gamma}^m\alpha_k)=b+pk$ for all $m\in \bZ$. Let
$\alpha_{k+me_K}=(-pu_{\gamma})^m\alpha_k$. So $v_L(\alpha_k)=b+pk$
for all $k\in \bZ$ and $\alpha_{k+e_K}=
-\sum_{i=1}^{p-1}\binom{p}{i}(\gamma-1)^{i-1}\alpha_k$.  As a result,
$(\gamma-1)\alpha_{k+e_K}= (\gamma-1)^p\alpha_k$, because
$(\gamma-1)^p=-\sum_{i=1}^{p-1}\binom{p}{i}(\gamma-1)^i$.

Now $v_L((\gamma-1)^i\alpha_k)=(i+1)b+pk$ for $0\leq i\leq
p-1$.  Use \cite[V \S3 Lem 4]{serre:local} to find $\mu_{i,k}\in N$
with $v_N(\mu_{i,k})=(1+pi)b+p^2k$ and
$\Phi_p(\sigma)\mu_{i,k}=(\gamma-1)^i\alpha_k$.  Since
$\Phi_p(\sigma)(\Theta-1)=(\gamma-1)\Phi_p(\sigma)$,
$\Phi_p(\sigma)\cdot\left ((\Theta-1)\mu_{i,k}- \mu_{i+1,k}\right )=0$
for $0\leq i\leq p-2$. Also $\Phi_p(\sigma)\cdot \left 
((\Theta-1)\mu_{p-1,k}-\mu_{1,k+e_K}\right )=0$.

By the Normal Basis Theorem, if $\Phi_p(\sigma)\nu=0$ for $\nu\in N$
with $v_N(\nu)\not\equiv b\bmod p$, then there is a $\theta\in N$ with
$v_N(\theta)=v_N(\nu)-b$ and $(\sigma-1)\theta=\nu$.  Recall
$u(\sigma)\in\euO_T[\sigma]^*$ defined in Lemma 13, and use the
Normal Basis Theorem to find $\rho_s^*\in N$ with $v_N(\rho_s^*)=b+ps$
such that
$$(\Theta-1)\mu_{i,k}=(\sigma-1)u(\sigma)\rho^*_{ib+pk+c}+\begin{cases}
\mu_{i+1,k}&\mbox{for }0\leq i<p-1\\
\mu_{1,k+e_K}&\mbox{for }i=p-1
\end{cases}$$
Now define $\rho_s\in N$ with $v_N(\rho_s)=b+ps$ by
$\rho_{bi+pk}=\mu_{i,k}$. And use Lemma 13 to replace 
$(\sigma-1)u(\sigma)\rho^*_s$ by 
$(1/p)\cdot ((\Theta-1)^p\rho_s^*+\sum_{i=1}^{p-1}\binom{p}{i}(\Theta-1)^i\rho_s^*)$.
\qed\end{pf}

\begin{pf*}{PROOF (Theorem 12 when $b_*<pb$)}
We have an $\euO_T$-basis for $\euP_N^r$ consisting of the elements
$$ p\rho_m, (\Theta-1) \rho_m, (\Theta-1)^2 p \rho_m, \ldots,
(\Theta-1)^{p-1}p\rho_m $$ 
for 
$$ \frac{r-b}{p} -pe_K \leq m < \frac{r}{p} + c -b_* \qquad
\mbox{(Range A)}; $$
and 
$$ (\Theta-1)^{j_m+1} \rho_m^*, \ldots, (\Theta-1)^{p-1} \rho_m^*,
p\rho_m , \ldots, (\Theta-1)^{j_m} p \rho_m $$
for 
$$ \frac{r}{p} + c -b_* \leq m < \frac{r-b}{p} \qquad \mbox{(Range
    B)}. $$
Here $m$ is restricted to integer values, and
$j_m \in \{0,\ldots, p-2\}$ is such that $(\Theta-1)^{j_m+1}\rho_m^*
\in \euP_N^r$ but $(\Theta-1)^{j_m}\rho_m^* \not \in \euP_N^r$.

Clearly $\mathcal{M}$ is spanned over $\bF_q$ by (the images of) these
$\euO_T$-basis elements.  So $\mathcal{M}$ is spanned over
$\bF_q[\Theta]$ by the $p\rho_m$ for $m \in$ Range A, along with the
$(\Theta-1)^{j_m+1}\rho_m^*$, and the $p \rho_m$ for $m \in$ Range B.
We want to determine the $\bF_q$-dimension of
$(\Theta-1)^{p-1}\mathcal{M}$.  So apply $(\Theta-1)^{p-1}$.  Since
the image of $(\Theta-1)^{p-1}p\rho_m$ in $\mathcal{M}$ is clearly
zero for $m\in$ Range B, we are left with an $\bF_q$-generating set
for $(\Theta-1)^{p-1}\mathcal{M}$ consisting of
\begin{equation}
(\Theta-1)^{p-1}p\rho_m\mbox{ for $m \in$ Range A, and }
(\Theta-1)^{j_m+p}\rho_m^*\mbox{ for $m \in$ Range B.}
\end{equation} 
This set of $e_K$ elements is not a basis.  Using the relationships in
Lemma 14, we should be able to replace certain
$(\Theta-1)^{p-1}p\rho_m$ with $(\Theta-1)^{j_m+p}\rho_{m'}^*$ for
some $m'\in$ Range B, or eliminate it entirely. Of course, since the
relationships in Lemma 14 are the only ``extra'' relations, once we
have made all such replacements/eliminations, we will be left with a
$\bF_q$-basis for $(\Theta-1)^{p-1}\mathcal{M}$.

Split Range B into a disjoint union of sets Range $\mbox{B}_0, \ldots
,$ Range $\mbox{B}_{p-2}$ where Range $\mbox{B}_j$ consists of those
$m$ with $j_m=j$.  In other words, $\mbox{Range }B_j=\{m\in\bZ:
r-(j+1)b_*-b\leq pm< r-jb_*-b\}$.

Take the relationship in Lemma 14, replace $m$ with $m-c$ and multiply
it by $(\Theta-1)^{p-2}$. We are interested in the situation where
$m-c\in$ Range A and $m\in$ Range B. Since the ``length'' of Range
B$_{p-2}$ is $b_*/p>c$, this actually occurs when $m-c\in$ Range A and
$m\in$ Range B$_{p-2}$. Since $(\Theta-1)^{i+p-2}p\rho_{m}^*\in
p\euP_N^r$ for $i\geq 1$ and $m\in$ Range B$_{p-2}$, the relationship
in Lemma 14 simplifies to
\begin{equation}
(\Theta-1)^{2p-2}\rho_{m}^* = 
(\Theta-1)^{p-1}p\rho_{m-c}-(\Theta-1)^{p-2}p\rho_{f(m)-c},
\end{equation}
where
\[f(m)=m+\begin{cases}b &\mbox {if }m\not\equiv c-b \bmod p,\\
pe_K-(p-2)b &\mbox {if }m\equiv c-b \bmod p.
\end{cases}\]

It is helpful, since we are interested in other relationships similar
to (4), to observe that in general,
\[v_N((\Theta-1)^{j+p}\rho_{m}^*) = v_N((\Theta-1)^{j+1}p\rho_{m-c})<v_N((\Theta-1)^jp\rho_{f(m)-c).}\]
So in regards to (4) where $j=p-2$, $(\Theta-1)^{p-2}p\rho_{f(m)-c}$
should be regarded as ``error.''  We can remove
$(\Theta-1)^{p-1}p\rho_{m-c}$ from the set of generators (3) for those
$m-c\in$ Range A such that $m\in$ Range B$_{p-2}$ and
$(\Theta-1)^{p-2}p\rho_{f(m)-c}\in p\euP_N^r$.  But before we do so,
need to consider $(\Theta-1)^{p-2}p\rho_{f(m)-c}\not\in p\euP_N^r$.
Indeed we will find that we do not need to treat these cases
separately.

Notice that $f(m)-m\geq b>b_*/p$, which is the approximate ``length''
of each Range B$_j$. So for $m\in$ Range B$_{p-2}$, it is certainly
the case that $f(m)\in$ Range B$_k$ for some $k\leq p-3$.  Moreover if
we denote iteration in the usual way, $f^2=f\circ f$, $f^3=f\circ
f\circ f$, etc., it will be the case that $f^2(m)\in$ Range B$_k$ for
some $k\leq p-4$ and so on. Think of Range B$_{-1}$ as including those $m$ such
that $p\rho_m^*\in p\euP_N^r$ and so are zero in $\mathcal{M}$.

Take the relationship in Lemma 14, replace $m$ with $m-c$ and multiply
by $(\Theta-1)^j$.  For $m\in$ Range B$_k$ with $k\leq j$, we have
$(\Theta-1)^{i+j}p\rho_m^*\in p\euP_N^r$ for $i\geq 1$. And so the
relationship simplifies to
$(\Theta-1)^{j+p}\rho_{m}^* = 
(\Theta-1)^{j+1}p\rho_{m-c}-(\Theta-1)^jp\rho_{f(m)-c}$.
As a result, in addition to (4), we also have
\begin{eqnarray*}
(\Theta-1)^{2p-3}\rho_{f(m)}^* &=& (\Theta-1)^{p-2}p\rho_{f(m)-c}-(\Theta-1)^{p-3}p\rho_{f^2(m)-c},\\
(\Theta-1)^{2p-4}\rho_{f^2(m)}^* &=& (\Theta-1)^{p-3}p\rho_{f^2(m)-c}-(\Theta-1)^{p-4}p\rho_{f^3(m)-c},\\
&\vdots &\\
(\Theta-1)^p\rho_{f^{p-2}(m)}^* &=& (\Theta-1)p\rho_{f^{p-2}(m)-c}-p\rho_{f^{p-1}(m)-c}.
\end{eqnarray*}
As a result, for $m-c\in$ Range A and $m\in$ Range B$_{p-2}$ we have
\[\sum_{j=0}^{p-2}(\Theta-1)^{p+j}\rho^*_{f^{p-j-2}(m)}= (\Theta-1)^{p-1}p\rho_{m-c}-p\rho_{f^{p-1}(m)-c},\]
where for $j<p-2$ either $f^{p-j-2}(m)\in$ Range B$_j$ and
$(\Theta-1)^{p+j}\rho^*_{f^{p-j-2}(m)}$ is a nontrivial generator
listed in (3), or $(\Theta-1)^{p+j}\rho^*_{f^{p-j-2}(m)}=0$ in
$\mathcal{M}$.  In any case, $(\Theta-1)^{p-1}p\rho_{m-c}$ is clearly
expressed in terms of other generators and can be removed {\em if and
only if} $p\rho_{f^{p-1}(m)-c}\in p\euP_N^r$. Notice that
\[f^{p-1}(m)=m+\begin{cases} (p-1)b &\mbox{if }m\equiv c\bmod p, \\
pe_K &\mbox{if }m\not\equiv c\bmod p.\end{cases}\]

As a result, we can remove those elements
$(\Theta-1)^{p-1}p\rho_{m-c}$ for $m\in$ Range B$_{p-2}$, namely
\[\frac{r}{p}+c-b_*\leq m< \frac{r+b}{p}+2c-b_*\]
such that $pf^{p-1}(m)\geq r+b_*-2b$. Once we have done so, we will have
an $\bF_q$-basis for $(\Theta-1)^{p-1}\mathcal{M}$.

Since $f^{p-1}(m)=m+pe_K$ for $m\not\equiv c\bmod p$ we can remove all
$m\not\equiv c\bmod p$. We can also remove all $m\equiv c\bmod p$ if
$(p^2-p+1)b\geq pb_*$. Doing so and keeping track of how many elements
were removed yields part of the statement of Theorem 12.  To get the
statement under near maximal refined ramification, notice that we need to
``put back'' one element for each integer $m\equiv c\bmod p$, that
satisfies $r/p+c-b_*\leq m$ and $pf^{p-1}(m)< r+b_*-2b$.
\qed\end{pf*}

We now state two corollaries of Theorem 12.
\begin{cor}
Let $K$ be a finite extension of $\mathbb{Q}_p$ and let $N_1, N_2$ be two fully ramified bicyclic extensions with unique ramification break number $b$.
Assume that the two second
refined ramification breaks satisfy
$b_*^{(1)}, b_*^{(2)}<(p-1+1/p)b$. If
$b_*^{(1)}\neq b_*^{(2)}$, then for each $r$,
$\euP_{N_1}^r\not\cong \euP_{N_2}^r$ as $\euO_T[G]$-modules.
\end{cor}
Motivated by the diagram in \S2.2.3, we observe that
when the break number $b$ is large enough, the hypothesis on the
second refined ramification numbers can be replace with a hypothesis
on $b$. 

\begin{cor}
Let $K$ be a finite extension of $\mathbb{Q}_p$ and let $N_1, N_2$ be two fully ramified bicyclic extensions with unique ramification break number $b$ satisfying
$$\left (1-\frac{p^2-2p+1}{p^3-2p+1}\right )\cdot \frac{pe_K}{p-1}<b<
\frac{pe_K}{p-1}.$$ If the two second
refined ramification breaks are different, $b_*^{(1)}\neq b_*^{(2)}$, then for each $r$,
$\euP_{N_1}^r\not\cong \euP_{N_2}^r$ as $\euO_T[G]$-modules.
\end{cor}

\subsection{Maximal refined ramification and Galois module structure}

In this section we assume $b_*=pb$ and establish an explicit integral
basis for $\euP_N^r$ over $\euO_T$ upon which we can follow the
Galois action in a particularly transparent way.

Recall the notation of \S2, in particular
$\Theta=\gamma\sigma^{[\omega]}$. Using Lemma 3, $b_*=pb$ means that
$v_N((\Theta-1)\rho= v_N(\rho)+pb$ for $\rho\in N$ if $v_N(\rho)\equiv
b\bmod p$ but $v_N(\rho)\not\equiv (1-p)b\bmod p^2$.  As a result, given
any $\rho_m\in N$ with $v_N(\rho_m)=b+p^2m$, we have
$v_N((\Theta-1)^i\rho_m)=(1+ip)b+p^2m$ for $0\leq i\leq p-1$ and thus
\begin{equation}
\rho_m^{(i, j)}:=(\Theta-1)^i(\sigma-1)^j\rho_m\mbox{ satisfies
}v_N(\rho_m^{(i, j)})=(1+j+ip)b+p^2m
\end{equation}
for $0\leq i,j \leq p-1$.
Since $\{v_N(\rho_m^{(i, j)}): 0\leq
i,j\leq p-1\}$ is a complete set of residues modulo $p^2$,
$\{v_N(\rho_m^{(i, j)}): 0\leq i,j\leq p-1, m\in \bZ\}=\bZ$ and since
$N/T$ is fully ramified, we can use the $\rho_m^{(i, j)}$ to construct
an $\euO_T$-basis for $\euP_N^r$. For example, simply choose
$\rho_m^{(i,j)}$ with $r\leq v_N(\rho_m^{(i, j)})\leq r+p^2e_K-1$.

So that the Galois action can be followed on this basis, we must
modify this construction (but only slightly).  Consider the `exponent'
$(i, j)$ to be a two digit $p$-ary integer $ip+j$.  The larger the
integer $(i, j)$ then, the larger the valuation $v_N(\rho_m^{(i,
j)})$.  Furthermore recall the diagram in \S2.2.3.  Since
$b_*=pb$ we have $b<p^2e_K/(p^2-1)$.  So there are values of $m$ such
that $r\leq v_N(p\rho_m^{(0, 0)})<v_N(p\rho_m^{(p-1,
p-1)})<r+p^2e_K$. For these $m$ set $(i, j)_m=(p-1, p-1)$.  Otherwise
$r\leq v_N(\rho_m^{(p-1, p-1)})< v_N(p\rho_m^{(0, 0)})<r+p^2e_K$.  For
each of these other values of $m$, let $(i, j)_m$ be the $p$-ary
integer such that $r\leq v_N(\rho_m^{(i, j)_m+(0, 1)})\leq
v_N(\rho_m^{(p-1, p-1)})< v_N(p\rho_m^{(0, 0)})\leq v_N(p\rho_m^{(i,
j)_m})<r+p^2e_K$. For each integer $m$ such that $(r-b)/p^2-e_K\leq
m<(r-b)/p^2$ define
\[\mathcal{M}(m)=
\euO_T\rho_m^{(i, j)_m+(0, 1)}+\cdots +\euO_T\rho_m^{(p-1,
p-1)}+\euO_T p\rho_m^{(0, 0)} +\cdots +\euO_Tp\rho_m^{(i, j)_m}.\]
Note that when $(i, j)_m=(p-1, p-1)$ we consider the sum
$\euO_T\rho_m^{(i, j)_m+(0, 1)}+\cdots +\euO_T\rho_m^{(p-1, p-1)}$ to
be empty.  In other words, $(p-1,p-1)+(0,1)$ should be considered larger
than $(p-1,p-1)$.  We find 
\[\euP_N^r=\sum_{m=A_r-e_K}^{A_r-1} \mathcal{M}(m)\mbox{ where }
A_r=\left\lceil \frac{r-b}{p^2}\right\rceil\] and $\lceil\cdot\rceil$
denotes the least integer function (ceiling function).

The $\euO_T[G]$-structure of $\euP_N^r$, namely Theorem 18, follows
then from the following lemma and some basic combinatorics.

\begin{lem} Each $\mathcal{M}(m)$ is 
isomorphic to an ideal of $\euO_T[G]$. Indeed, if we
write $(i, j)_m$ as $(i_m, j_m)$, then
$$\mathcal{M}(m)\cong
\left\langle p,
(\Theta-1)^{i_m}(\sigma-1)^{j_m+1},
(\Theta-1)^{i_m+1}\right\rangle$$
\end{lem}
\begin{pf}
Let $\phi_{(i, j)}=(\Theta-1)^i(\sigma-1)^j$. So $\rho_m^{(i,
j)}=\phi_{(i, j)}\rho_m$. Now list the $\euO_T$-basis elements of
$\mathcal{M}(m)$ in an array, dropping the $\rho_m$ from each element:
\begin{equation}
\begin{array}{c}
\begin{array}{ccc|ccc}
p\phi_{(0, 0)}, & \cdots & p\phi_{(i_m, 0)}, & \phi_{(i_m+1, 0)}, &\cdots &
 \phi_{(p-1, 0)},\\ & \ddots & & & \ddots & \\ 
p\phi_{(0, j_m)}, &\cdots & \hspace*{5.5mm}p\phi_{(i_m, j_m)},\hspace*{6mm} & \phi_{(i_m+1, j_m)}, &\cdots
 & \phi_{(p-1, j_m)}, \\ \cline{3--4}
\end{array} \\
\begin{array}{cc|cccc}
\hspace*{4mm}p\phi_{(0, j_m+1)}, & \cdots &
\hspace*{2mm}\phi_{(i_m, j_m+1)}, & \hspace*{4mm}\phi_{(i_m+1, j_m+1)}, &\cdots &
\hspace*{4mm}\phi_{(p-1, j_m+1)},\\ & \ddots & &  &\ddots & \\
p\phi_{(0, p-1)}, & \cdots &
\hspace*{2mm}\phi_{(i_m, p-1)}, & \phi_{(i_m+1, p-1)}, & \cdots &
\phi_{(p-1, p-1)}.
\end{array} 
\end{array}
\end{equation}
The boundary between elements in and out of $p\euO_T[G]$ is marked.
We would like to show that the $\euO_T$-span of (6) is the
$\euO_T[G]$-ideal generated by $p\phi_{(0, 0)}$, $\phi_{(i_m,
j_m+1)}$, $\phi_{(i_m+1, 0)}$. But this is clear, once we know that
this $\euO_T$-span is closed under $\sigma$ and $\Theta$, and this
follows from the fact that $(\sigma-1)^p$ and
$(\Theta-1)^p\in p\euO_T[G]$.  \qed\end{pf}

\begin{thm} 
Let $b_*=pb$.  For $0\leq s\leq p^2-1$ let $(i_s, j_s)$ denote the
$p$-ary expansion of $s$. So $s=i_sp+j_s$.  For $0\leq s\leq p^2-1$,
let $\mathcal{I}_s$ be the $\euO_T[G]$ ideal $\left\langle p,
(\Theta-1)^{i_s}(\sigma-1)^{j_s+1}, (\Theta-1)^{i_s+1}\right\rangle$, with 
$\mathcal{I}_{p^2-1}= p\euO_T[G]\cong \euO_T[G]$.  Then
\[\euP_N^r\cong 
\bigoplus_{s=0}^{p^2-2} \mathcal{I}_s
^{\left\lceil\frac{r-(s+1)b}{p^2}\right\rceil-
\left\lceil\frac{r-(s+2)b}{p^2}\right\rceil}
\oplus
\mathcal{I}_{p^2-1}^{e_K-b+
\left\lceil\frac{r}{p^2}\right\rceil-\left\lceil\frac{r-b}{p^2}\right\rceil}
\] 
as
$\euO_T[G]$-modules. This structure is
parametrized by $r, \omega, p, e_K$ and $b$.
\end{thm}

\begin{pf*}{PROOF (Theorem 12 when $b_*=pb$).}
Each $\euO_T[G]$-module $\mathcal{I}_s$ has an $\euO_T$-basis as in
(6).  To compute
$\mbox{dim}_{\mathbb{F}_q}((\Theta-1)^{p-1}\mathcal{I}_s/p\mathcal{I}_s)$,
simply apply $(\Theta-1)^{p-1}$ to each of these $\euO_T$-basis
elements.  By Lemma 13, for $i\geq 1$ we have
$(\Theta-1)^{p-1}\cdot(\Theta-1)^i(\sigma-1)^j\equiv u(\sigma)\cdot
(\Theta-1)^{i-1}(\sigma-1)^{j+1}p \bmod p(\Theta-1)^i(\sigma-1)^j$,
where $u(\sigma)$ is a unit.  As a result,
$$\mbox{dim}_{\mathbb{F}_q}((\Theta-1)^{p-1}\mathcal{I}_s/p\mathcal{I}_s)=
\begin{cases}
p-1&0\leq s\leq p-1,\\
p-2&p\leq s\leq p^2-p-1, s\not\equiv -1\bmod p,\\
p-1&p\leq s\leq p^2-p-1, s\equiv -1\bmod p,\\
p-1&p^2-p\leq s\leq p^2-2,\\
p&s= p^2-1.
\end{cases}$$
Let $m_s$ denote the multiplicity of $\mathcal{I}_s$ in the statement
of Theorem 18.  So $m_s=\lceil (r-(s+1)b)/p^2\rceil-\lceil
(r-(s+2)b)/p^2\rceil$ for $0\leq s\leq p^2-2$. Therefore because
$p-1=(p-2)+1$, we have
$\mbox{dim}_{\mathbb{F}_q}((\Theta-1)^{p-1}\euP_N^r/p\euP_N^r)=pm_{p^2-1}+
(p-2)\sum_{s=0}^{p^2-2}m_s+\sum_{s=0}^{p-1}m_s+\sum_{s=p^2-p}^{p^2-2}m_s+
\sum_{k=2}^{p-1}m_{kp-1}$. These are for the most part telescoping
sums, and so the expression simplifies to
$\mbox{dim}_{\mathbb{F}_q}((\Theta-1)^{p-1}\euP_N^r/p\euP_N^r)=pe_k-2b+
\lceil r/p^2\rceil-\lceil (r-b)/p^2\rceil -\lceil (r-(p+1)b)/p^2\rceil
+\lceil (r+(p-1)b)/p^2\rceil+\sum_{k=2}^{p-1}m_{kp-1}$. It remains to
recognize that $\sum_{k=0}^{p-1}m_{kp-1}=\lceil
r/p\rceil-\lceil(r-b)/p\rceil$, which follows from the fact that both
count the number of integers $i$ such that $(r-b)/p\leq i\leq
(r-1)/p$.  Each term $m_{kp-1}$ in the sum simply counts those
integers $\equiv kb\bmod p$.
\end{pf*}

\section{Conclusion}

This paper grows out of the on-going effort to generalize the
biquadratic results of \cite{elder:onebreak} to $p>2$. Thus far,
several themes have emerged
and a number of questions have been raised, all of which which bear repeating.

The central theme is the role of truncated exponentiation. Its
appearance within the group ring $\euO_T[G]$ for $G$ elementary
abelian has led to the refined ramification filtration
\cite{elder:newbreaks}. Notably, the definition of refined
ramification break numbers remains tied to a choice of element and so
the refined breaks (beyond the first two) cannot, as yet, be said to
be canonical. In addition, it has been observed in the context of
quaternion extensions \cite{elder:hooper} that the refined
ramification filtration has some influence on breaks in the usual
ramification filtration. And so there is much remaining work to determine
if/how these two filtrations fit together.

The appearance of truncated exponentiation among the generators of the
extension in \S2.2.2 as well as the notion of maximal refined
ramification and the ease and transparency in \S3.2 are the motivation
for \cite{elder:onedim}. This, along any connection with Artin-Hasse
exponentiation and explicit reciprocity \cite{fesenko}, warrant
further investigation.

In \cite[\S4]{elder:onebreak} and then here in \S2.2.3, a question is
raised concerning how twists by characters of Galois representations
effect ramification, refined ramification and Galois module structure.
One consequence of this question is the suggestion that the problem of
Galois module structure be broken in two: (1) The determination of
nice classes of extension, for which the Galois module structure can
be easily determined. {\em e.g.} \cite{elder:onedim} (2) The problem
of Galois module structure under twisting, which remains very much
open.

\bibliography{bib}





\end{document}